\documentstyle[11pt]{article}
\newlength{\standardunitlength}
 \newtheorem{lemma}{Lemma}
\newtheorem{theorem}{Theorem} \newtheorem{prop}{Proposition}
\newenvironment{proof}{\noindent {\sc Proof:}}{$\Box$ \vspace{2 ex}}

\begin{document}

\begin{center}
Descent Identities, Hessenberg Varieties, and the Weil Conjectures
\end{center}

\begin{center}
By Jason Fulman
\end{center}

\begin{center}
Dartmouth College
\end{center}

\begin{center}
Jason.E.Fulman@Dartmouth.Edu
\end{center}

\begin{center}
Abstract: The Weil Conjectures are applied to the Hessenberg Varieties to obtain
interesting information about the combinatorics of descents in the
symmetric group. Combining this with elementary linear algebra leads to
elegant proofs of some identities from the theory of descents.
\end{center}

\section{Introduction and Background}

	The purpose of this introduction is to give background on the
following three topics: permutation statistics, Hessenberg varieties, and
the Weil conjectures. The topics will be described in this order, and the
emphasis will be on explaining their relationship to each other as is
relevant to this note.

	Permutation statistics are functions from the symmetric group $S_n$
to the non-negative integers. Many permutations statistics have interesting
combinatorial properties (pages 17-31 of Stanley \cite{Stanley}) and give
rise to metrics which are important in the statistical theory of ranking
(Chapter 6 of Diaconis \cite{Diaconis}). Volume 3 of Knuth \cite{Knuth}
connects permutation statistics with the theory of sorting.

	One important statistic is the number of inversions of a
permutation. This is denoted $Inv(\pi)$ and is equal to the number of pairs
$(i,j)$ such that $1\leq i < j \leq n$ and $\pi(i)>\pi(j)$. The number of
inversions of $\pi$ is also equal to the length of $\pi$ in terms of the
generating reflections $\{(i,i+1):1 \leq i \leq n-1 \}$. Inversions have
the following well-known generating function (e.g. page 21 of Stanley
\cite{Stanley})

\[ \sum_{\pi \in S_n} q^{Inv(\pi)} = \prod_{i=1}^n \frac{q^i-1}{q-1} \]

which can be used to prove that the distribution of inversions is
asymptotically normal in the $n \rightarrow \infty$ limit (e.g. Chapter 6
of Diaconis \cite{Diaconis}). It is worth observing that $\prod_{i=1}^n
\frac{q^i-1}{q-1}$ is equal to the number of complete flags (i.e. $V_0=id
\subset V_1 \subset \cdots \subset V_n=V$ with $dim(V_i)=i$) in an
$n$-dimensional vector space over a finite field of size $q$. This already
suggests a connection with algebraic geometry.

	Inversions can be defined for any finite reflection group $W$. For
$w \in W$, define $Inv(w)$ as the length of $w$ in terms of the simple
reflections (i.e. those corresponding to simple roots). Alternately,
$Inv(w)$ is the number of positive roots which $w$ send to negative
roots. Here too there is a factorization

\[ \sum_{w \in W} q^{Inv(w)} = \prod_{i=1}^n \frac{q^{d_i}-1}{q-1} \]

where $d_i$ are the degrees of the corresponding reflection group (see page
73 of Humphreys \cite{Humphreys}).

	A second permutation statistic of interest is the number of
descents of a permutation. This is denoted $d(\pi)$ and is equal to the
number of pairs $(i,i+1)$ with $1 \leq i \leq n-1$ and $\pi(i) >
\pi(i+1)$. The generating function for descents gives rise to the Eulerian
polynomials

\[ A_n(q) = \sum_{\pi \in S_n} q^{d(\pi)+1} \]

	Pages 243-246 of Comtet \cite{Comtet} describe some properties of
the Eulerian polynomials. Two of the nicest generating functions involving
Eulerian polynomials are

\[ \frac{A_n(q)}{(1-q)^{n+1}} = \sum_{l \geq 0} l^nq^l \]

\[ \sum_{n \geq 0} \frac{A_n(q) t^n}{n!} = \frac{1-q}{1-qe^{t(1-q)}} \]

	Some recurrences for the Eulerian polynomials will be found in
Section $\ref{Descent Identities}$ of this paper. Descents can be defined
for all Coxeter groups as the number of simple positive roots $w$ makes
negative. In general no nice generating function seems to be known.

	Next, we recall the Hessenberg varieties defined by DeMari,
Procesi, and Shayman \cite{DeMari}. Let $G$ be a complex, semisimple
algebraic group, $B$ a Borel subgroup, and $T$ a maximal torus in $B$. Let
$\bf{g,b,t}$ be the lie algebras of $G,B,T$ respectively. Let $\bf{h}$ be a
subspace of $\bf{g}$ which contains $\bf{b}$ and is a $\bf{b}$
submodule. Let $s \in \bf{t}$ be a regular, semi-simple element. Define the
corresponding Hessenberg variety (which turns out to depend on $G$ and
$\bf{h}$ but not on $s$) by

	\[ X_H(s) = \{ g \in G | Ad(g^{-1})[s] \in \bf{h} \}\]

where $Ad$ is the Adjoint action of Lie theory (conjugation in the case of
matrix groups).

	The main example to be considered in this note is the following,
which we will call $Hess(n,p)$. Fix $p: 1 \leq p \leq n-1$. Let $G=SL(n,C)$
and $\bf{h}$ be the subspace of $sl(n,C)$ consisting of those matrices
$(h_{ij})$ for which $h_{ij}=0$ if $i-j>p$. Let $s$ be any diagonalizable
element of $G$ with distinct eigenvalues. Then the corresponding Hessenberg
variety $X_H(s)$ can be more simply described as all complete flags $V_0
\subset V_1 \subset \cdots \subset V_n$ satisfying the condition that
$s(V_i) \subset V_{i+p}$. For example, $Hess(n,n-1)$ is the flag variety of
$SL(n,C)$.

	DeMari, Procesi, and Shayman \cite{DeMari} study the varieties
$X_H(s)$, proving that they are smooth toric varieties and computing their
Betti numbers. We will require only the following special case. The point
for us is that the Betti numbers of the varieties $Hess(n,p)$ are
interesting permutation statistics.

\begin{theorem} (DeMari,Shayman \cite{hess}) \label{bet}

\begin{enumerate}
\item The varieties $Hess(n,p)$ are smooth.
\item The odd Betti
numbers $\beta_{2k-1}(Hess(n,p))$ vanish. The even Betti numbers
$\beta_{2k}(Hess(n,p))$ are equal to the number of permutations on $n$
symbols such that $|\{(i,j):1 \leq i < j \leq n, j-i \leq p,
\pi(i)>\pi(j)\}|=k$.
\item For $q$ sufficiently large, the equations defining $Hess(n,p)$,
reduced to a field of $q$ elements, define a smooth variety.
\end{enumerate}

\end{theorem}

	The third part of Theorem \ref{bet} was not explicitly stated by
DeMari and Shayman \cite{hess}, but follows by the same arguments as in the
smooth case, given on pages 224-5 of their article. If $q \leq n$ then
there does not exist an invertible $n*n$ diagonal matrix with distinct
eigenvalues all of which lie in a field of $q$ elements. Throughout the
rest of this paper it will be assumed that $q$ is sufficiently large so that the
reduced variety is smooth.

	As an example of Theorem \ref{bet}, $\beta_{2k}(Hess(n,n-1))$ is
the number of permutations in $S_n$ with $k$ inversions and
$\beta_{2k}(Hess(n,1))$ is the number of permutations in $S_n$ with $k$
descents.

	Next let us review the Weil conjectures. One use of them is to
compute Betti numbers of continuous varieties by counting points in
varieties defined over finite fields. This will be done in Section
\ref{Descent Identities}, thereby proving identities about the Eulerian
polynomials. The version of the Weil conjectures considered here can be
found in Appendix C of Hartsthorne \cite{Hart}. These conjectures are now,
of course, theorems.

\begin{theorem} \label{Weil} (Weil Conjecture) Given a smooth variety $V$,
its Betti numbers can be computed as follows. Reduce the equations defining
$V$ to equations over a field of $q^s$ elements where $q$ is a prime power,
and let $N(q^s)$ be the number of solutions to these reduced equations. If
the reduced variety is smooth for all such reductions then

	\[ exp \sum_{s=1}^{\infty} \frac{N(q^s)x^s}{s} =
\frac{P_1(x)P_3(x)...P_{2\delta-1}(x)}{P_0(x)P_2(x)...P_{2\delta}(x)} \]

where $\delta$ is the dimension of $V$ and $P_k(x)=\prod_{j=1}^{\beta_k}
(1-\alpha_{kl}x)$, with $|\alpha_{kl}|=q^{\frac{k}{2}}$.
\end{theorem}

	Stanley \cite{StanleySchub} has written Theorem \ref{Weil} in a
form which is somewhat more useful for our purposes.

\begin{prop} (Stanley \cite{StanleySchub}) \label {simplify} Suppose in
addition to the assumptions of Theorem \ref{Weil} that $N(q^s)$ is a
polynomial $\sum_k \gamma_k q^{ks}$ in $q^s$. Then $\beta_{2k}=\gamma_k$.
\end{prop}

\begin{proof}
\begin{eqnarray*}
exp \sum_{s=1}^{\infty} \frac {N(q^s)x^s}{s} & = & exp \sum_{s=1}^{\infty} (\sum_k
\gamma_k q^{ks}) \frac{x^s}{s}\\
& = & exp \sum_k \gamma_k \sum_{s=1}^{\infty} \frac{(q^kx)^s}{s}\\
& = & exp \sum_k -\gamma_k ln(1-q^kx)\\
& = & \prod_k (1-q^kx)^{-\gamma_k}
\end{eqnarray*}
\end{proof}

\section{Descent Identities} \label{Descent Identities}

	As an example of the concepts in the introduction, let us use the
Weil conjectures to find the generating function for permutations in $S_n$
by inversions (also known as the Poincar\'e series of $S_n$).

\begin{theorem} \label{inv}
\[ \sum_{\pi \in S_n} q^{Inv(\pi)} = \prod_{i=1}^n \frac{q^i-1}{q-1} \]
\end{theorem}

\begin{proof}
	Theorem $\ref{bet}$ and Proposition \ref{simplify} applied to
$Hess(n,n-1)$ show that

\[ \sum_{\pi \in S_n} q^{Inv(\pi)} = N(q),\]

the number of complete flags $V_0 \subset V_1 \subset ... \subset
V_n$ over a field of $q$ elements. $V_1$ can be chosen in
$\frac{q^n-1}{q-1}$ ways. Given this choice of $V_1$, modding out the
flag by $V_1$ shows that $V_2$ can be chosen in
$\frac{q^{n-1}-1}{q-1}$ ways. Continuing in this way and multiplying
proves that

\[ N(q)= \prod_{i=1}^n \frac{q^i-1}{q-1} \]

for infinitely many $q$, hence for all $q$ since both sides are
polynomials. 
\end{proof} 

	The above proof of Theorem \ref{inv} seems not to have been written
down before, perhaps because of the immense difficulty in proving the Weil
conjectures. Chevalley \cite{Chevalley} and Bott \cite{Bott} used the
topology of compact Lie groups to prove the factorization of the Poincar\'e
series for Weyl groups. The argument in Theorem \ref{inv} extends to the
other Weyl groups, but this would be somewhat circular because one must
know the size of the flag variety $G/B$ where $G$ is a finite algebraic
group with Weyl group $W$, and historically $|G|$ was computed using the
Bruhat decomposition and the factorization of the Poincar\'e series for
Weyl groups.

	The linear algebra involved in using the Weil conjectures to study
the Eulerian polynomials $A_n(q)$ is slightly more involved. We thus
establish two easy lemmas.

\begin{lemma} \label{subspace} Let $M \in GL(n,K)$ act on an $n$
dimensional vector space $V$ over a field $K$ such that $M$ has distinct
eigenvalues which are all contained in $K$. Then there are exactly ${n
\choose m}$ subspaces of dimension $m$ which are invariant under $M$.
\end{lemma}

\begin{proof}
	Let $W$ be a subspace of dimension $m$ which is invariant under
$M$. The characteristic polynomial of $M$ restricted to $W$ divides the
characteristic polynomial of $M$ on $V$, since $W$ is invariant. Since $M$
has distinct eigenvalues on $V$, its characteristic polynomial consists of
distinct linear factors, so the same is true for the characteristic
polynomial of $M$ on $W$. Thus $W$ is spanned by some $m$ of the $n$
1-dimensional eigenspaces for the action of $M$ on $V$.  \end{proof}

	Given a linear transformation $M$ on a $n$ dimensional vector space
$V$, call a vector $\vec{v}$ primitive if the set
$\{\vec{v},M\vec{v},M^2\vec{v},...,M^{n-1}\vec{v}\}$ forms a basis of $V$.

\begin {lemma} \label{primit} Let $M \in GL(n,K)$ act on an $n$ dimensional
vector space $V$ over a field $K$ such that $M$ has distinct eigenvalues
which are all contained in $K$. Then a vector $\vec{v}$ is primitive if and
only if its components with respect to a basis of eigenvectors of $V$ are
all non-0.
\end {lemma}

\begin{proof}
	Pick a basis of eigenvectors $\vec{e_1},...,\vec{e_n}$ of $M$ with
eigenvalues $\lambda_1,...,\lambda_n$. Let $\vec{v}$ have components
$(v_1,\cdots,v_n)$ with respect to this basis. Then $M^i \vec{v}$ has
components $(\lambda_1^iv_1,...,\lambda_n^iv_n)$. Clearly $\vec{v}$ is
primitive if and only if the determinant of the matrix with rows
$\vec{v},M\vec{v},M^2\vec{v},...,M^{n-1}\vec{v}$ written with respect to
the basis of eigenvectors, is non-vanishing. The value of this determinant
is

\[ [\prod_{i=1}^n v_i] det \left( \begin{array}{c c c c}
		1 & 1 & \cdots & 1 \\
		\lambda_1 & \lambda_2 & \cdots & \lambda_n \\
		\cdots & \cdots & \cdots & \cdots \\
		\lambda_1^{n-1} & \lambda_2^{n-1} & \cdots & \lambda_n^{n-1}
	  \end{array} \right) \]
 
	which is non-vanishing precisely when all $v_i$ are non-vanishing
by the theory of Vandermonde determinants and the fact that the eigenvalues
$\lambda_i$ of $M$ are distinct.
\end{proof}

	Recall that $A_n(q)$ denotes the $nth$ Eulerian polynomial
$\sum_{\pi \in S_n} q^{d(\pi)+1}$. For convenience set $A_0(q)=q$. Theorem
\ref{maybenew} is likely known, though we have not seen it in the
literature before.

\begin{theorem} \label{maybenew} $A_n(q)= \sum_{i=1}^n {n \choose
i}(q-1)^{i-1} A_{n-i}(q)$.
\end{theorem}

\begin{proof}
	Let $N(q,n)$ be the number of flags $0 \subset V_1 \subset V_2
\subset ... \subset V_n=V$ over the field of $q$ elements such that $MV_j
\subset V_{j+1}$ for all $j$, where $M \in GL(n,q)$ is a diagonal matrix
with distinct eigenvalues. Let $i$ be the smallest number between 1 and $n$
such that $MV_i=V_i$. By Lemma \ref{subspace}, there are ${n \choose i}$
ways of choosing $V_i$. The part of the flag between $V_1$ and $V_i$ is
determined by $V_1$, which is spanned by a primitive vector in the $i$
dimensional space $V_i$. There are, by Lemma \ref {primit}, $(q-1)^i$
primitive vectors for $V_i$, and hence $(q-1)^{i-1}$ choices for
$V_1$. Modding out the part of the flag between $V_i$ and $V_n=V$ by $V_i$
shows that there are $N(q,n-i)$ possibilities for this part of the flag. We
thus have the recurrence

\[ N(q,n)=\sum_{i=1}^n {n \choose i}(q-1)^{i-1}N(q,n-i) \]

	By Proposition \ref{simplify}, $A_n(q)=qN(q,n)$, proving the
theorem.
\end{proof}

	The recurrence in Theorem \ref{maybenew} was proved by splitting
the flag at the first subspace invariant under $M$ and summing over such
splittings. Theorem \ref{Frob} will come from splitting the flag at all the
$W_i$ invariant under $M$, and then summing over all such splittings. The
recurrence in Theorem \ref{Frob} is known and goes back to Frobenius
\cite{Frobenius}, though the proof is completely different. The notation
$S(n,r)$ means a Stirling number of the second kind and is the number of
set partitions of $\{1,\cdots,n\}$ into $r$ parts (page 33 of Stanley
\cite{Stanley}).

\begin{theorem} \label{Frob} $A_n(q)= q \sum_{r=1}^n r! S(n,r) (q-1)^{n-r}$.
\end{theorem}

\begin{proof}
	Proposition \ref{simplify} shows that $A_n(q)=qN(q,n)$, where
$N(q,n)$ is the number of flags $V_0=0 \subset V_1 \subset V_2 \subset
... \subset V_n=V$ such that $MV_j \subset V_{j+1}$ for $1 \leq j \leq n-1$
and $M \in GL(n,q)$ is a diagonal matrix with distinct eigenvalues. We
count these flags by the set $I$ of $i>0$ such that $V_i$ is invariant
under $M$. For each subset $I=\{i_1,i_2,...,i_r=n\}$ of $\{1,...,n\}$,
there are, by Lemma \ref{subspace}, ${n \choose i_{r-1}}{i_{r-1} \choose
i_{r-2}}...{i_2 \choose i_1}$ ways of picking the invariant subspaces
$V_{i_1},V_{i_2},\cdots,V_{i_r=n}$ of dimensions $\{i_1,i_2,...,i_r=n\}$ so
as to respect the inclusion relations. Consider the portion of the flag
between two consecutive invariant subspaces, say $V_{i_1} \subset V_{i_1+1}
\subset ... \subset V_{i_2}$. Modding out this whole sequence by $V_{i_1}$
shows that $V_{i_1+1}/V_{i_1}$ must be spanned by a primitive vector for
the action of $M$ on $V_{i_2}/V_{i_1}$. The dimension of the quotient is
$i_2-i_1$ so by Lemma \ref{primit} there are $(q-1)^{i_2-i_1}$ primitive
vectors. Since we are only interested in the subspace spanned by the
vector, we divide out by $q-1$. Multiplying out these choices of primitive
vectors, one sees that there are $(q-1)^{n-r}$ such choices.

	Therefore,

\begin{eqnarray*}
A_n(q) & = & q N(q,n)\\
& = & q \sum_{r=1}^n (q-1)^{n-r} \sum_{I \subset \{1,\cdots,n\} \atop n \in
I, |I|=r} {n \choose i_{r-1}}{i_{r-1} \choose i_{r-2}}...{i_2 \choose
i_1}\\
& = & q \sum_{r=1}^n r! S(n,r) (q-1)^{n-r}
\end{eqnarray*}

	The last equality follows because

\[ \sum_{I \subset \{1,\cdots,n\} \atop n \in I, |I|=r} {n \choose
n-i_{r-1}}{i_{r-1} \choose i_{r-1}-i_{r-2}}...{i_2 \choose i_2-i_1} \]

	is the number of ways of choosing a set partition of
$\{1,\cdots,n\}$ into $r$ parts with an ordering on the parts (the first
part has size $n-i_{r-1}$ and can be chosen in ${n \choose n-i_{r-1}}$
ways, the second part has size $i_{r-1}-i_{r-2}$ and can be chosen in
${i_{r-1} \choose i_{r-1}-i_{r-2}}$ ways, etc.) However by the definition
of the Stirling numbers of the second kind, the number of set partitions of
$\{1,\cdots,n\}$ into $r$ parts with an ordering on the parts is equal to
$r! S(n,r)$.
\end{proof}

\section{Concluding Thoughts}

	Here are some concluding thoughts about the results of this paper.

\begin{enumerate}

\item The proofs of the descent identities given here admittedly use a lot
of machinery. Nevertheless, given this machinery, the method of counting
flags employed in Theorems \ref{inv}, \ref{maybenew}, and \ref{Frob} is
natural and gives one a feel for where the recurrences come from. Direct
combinatorial proofs of these theorems would require more imagination.

	However, suppose one wants a recurrence for the Eulerian numbers
$A(n,k)$ which are the number of permutations on $n$ symbols with $k+1$
descents. It is easy to see combinatorially that
	
	\[ A(n,k) = (n-k+1)A(n-1,k-1) + kA(n-1,k) \]

	Thus direct combinatorics seems superior for finding recursions
satisfied by the coefficients of the Eulerian polynomials, but the flag
counting methods seem well-adapted toward finding recurrences satisfied by
the polynomials themselves.

\item It would be interesting to use the Weil conjectures to find
recurrences for the generating functions for the permutation statistics
whose value at $\pi$ is $|\{(i,j):1 \leq i < j \leq n, j-i \leq p,
\pi(i)>\pi(j)\}|$. Descents and inversions are the limiting cases $p=1,n-1$
and are the Betti numbers of $Hess(n,1)$ and $Hess(n,n-1)$ respectively. We
have not made much progress for other $p$ (though it is worth pointing out
that direct combinatorial arguments have not been successfully applied to
this problem either). It might also be interesting to study the generating
functions for descents for other Weyl groups.

\end{enumerate}

\bigbreak \section*{Acknowledgments} This work was done under the financial
support of an Alfred P. Sloan Dissertation Fellowship.

\end{document}